\documentclass[10pt,a4paper]{article}
\usepackage{cite}
\usepackage{url}
\usepackage{graphicx}
\usepackage{a4wide}
\usepackage[english]{babel}
\usepackage{amsfonts}
\usepackage{amsmath}
\usepackage{amssymb}
\usepackage{enumitem}
\usepackage{float}
\usepackage[latin1]{inputenc}
\usepackage{url}
\usepackage{tocbibind, times, theorem}
\usepackage{t1enc}
\usepackage{geometry}
\allowdisplaybreaks
\newtheorem{lemma}{Lemma}

\newtheorem{cor}{Corollary}
\newtheorem{theorem}{Theorem}

\DeclareMathOperator{\E}{\mathbb{E}}

\newcommand{\ve}{\varepsilon}
\makeatletter\def\blfootnote{\xdef\@thefnmark{}\@footnotetext}\makeatother

\title{\bf Probabilistic discrepancy bound for Monte Carlo point sets}
\author{Christoph Aistleitner\footnote{Graz University of Technology,
Institute of Mathematics A, Steyrergasse 30, 8010 Graz, Austria. \mbox{e-mail}: \texttt{aistleitner@math.tugraz.at}. Research supported by the Austrian Research Foundation (FWF),
Project S9603-N23.} \quad Markus Hofer\footnote{Graz University of Technology,
Institute of Mathematics A, Steyrergasse 30, 8010 Graz, Austria. \mbox{e-mail}: \texttt{markus.hofer@tugraz.at}. Research supported by the Austrian Research Foundation (FWF),
Project S9603-N23 and by the DOC [Doctoral Fellowship Programme of the Austrian Academy of Sciences].}}
\begin{document}

\date{}
\maketitle


\begin{abstract}
By a profound result of Heinrich, Novak, Wasilkowski, and Wo{\'z}niakowski the inverse of the star-discrepancy $n^*(s,\ve)$ satisfies the upper bound $n^*(s,\ve) \leq c_{\mathrm{abs}} s \ve^{-2}$. This is equivalent to the fact that for any $N$ and $s$ there exists a set of $N$ points in $[0,1]^s$ whose star-discrepancy is bounded by $c_{\mathrm{abs}} s^{1/2} N^{-1/2}$. The proof is based on the observation that a random point set satisfies the desired discrepancy bound with positive probability. In the present paper we prove an applied version of this result, making it applicable for computational purposes: for any given number $q \in (0,1)$ there exists an (explicitly stated) number $c(q)$ such that the star-discrepancy of a random set of $N$ points in $[0,1]^s$ is bounded by $c(q) s^{1/2} N^{-1/2}$ with probability at least $q$, uniformly in $N$ and $s$.
\end{abstract}

\section{Introduction and statement of results}

The number $n^*(s,\ve)$, which is defined as the smallest possible cardinality of a point set in $[0,1]^s$ having discrepancy bounded by $\ve$, is called the \emph{inverse of the discrepancy}. Heinrich, Novak, Wasilkowski, and Wo{\'z}niakowski \cite{hnww} proved the upper bound
\begin{equation} \label{inv}
n^*(s,\ve) \leq c_{\mathrm{abs}} s \ve^{-2},
\end{equation}
which is complemented by the lower bound
$$
n^*(s,\ve) \geq c_{\mathrm{abs}} s \ve^{-1}
$$
due to Hinrichs \cite{hc} (throughout the paper, $c_{\mathrm{abs}}$ denotes absolute constants, not always the same). Hence the inverse of the star-discrepancy depends linearly on the dimension, while the precise dependence on $\ve$ is still unknown. It is easy to see that \eqref{inv} is equivalent to the fact that for any $N$ and $s$ there exists a set $\mathcal{P}_N$ of $N$ points in $[0,1]^s$ such that the star-discrepancy $D_N^*$ of this point set is bounded by 
\begin{equation} \label{inv2}
D_N^* (\mathcal{P}_N) \leq c_{\mathrm{abs}} \frac{\sqrt{s}}{\sqrt{N}}
\end{equation}
(recently we showed that it is possible to choose $c_{\mathrm{abs}}=10$ in \eqref{inv2}, see \cite{ady}). The existence of such a point set directly follows from the surprising observation that a randomly generated point set (that is, a Monte Carlo point set) satisfies the desired discrepancy estimate with positive probability. Of course, for applications such a mere existence result is not of much use, as was remarked by several colleagues at the MCQMC 2012 conference in Sydney. For this reason, in the present paper we prove an applied version of \eqref{inv2}, which provides estimates for the probability of a random point set satisfying \eqref{inv2} (depending on the value of the constant). As our Theorem \ref{th1} below shows, this probability is extremely large already for moderate values of $c$, for example for $c=20$. Additionally, the quality of our estimates for these probabilities \emph{improves} as the dimension $s$ increases (which is somewhat counter-intuitive, and originates from the exponential inequalities used in the proof, which cause a ``concentration of mass'' phenomenon).\\

The fact that the probability of a random point set satisfying \eqref{inv2} is very large is in contrast to the fact that no general constructions of point sets satisfying such discrepancy bounds are known. So far, the best results are a component-by-component construction of Doerr, Gnewuch, Kritzer and Pillichshammer \cite{dgk}, a semi-deterministic algorithm based on dependent randomized rounding due to Doerr, Gnewuch, and Wahlstr{\"o}m \cite{dgw}, and a construction of Hinrichs \cite{hd} of a ``structured'' set of $N=256$ points in dimension $s=15$ having discrepancy less than 1/4  (by this means solving one instance of an open problem in \cite{now2}).\\

For more information concerning the inverse of the discrepancy and tractability of multidimensional integration we refer to a recent survey article of Gnewuch \cite{gnesurvey}, and to the monograph of Novak and Wo{\'z}niakowski \cite{now1, now2}. A collection of open problems on this topic can be found in \cite{hopen}.\\

In the present paper, we will prove the following theorem.

\begin{theorem} \label{th1}
For any $s \geq 1, ~N \geq 1$ and $q \in (0,1)$ a randomly generated $s$-dimensional point set $(z_1, \ldots, z_N)$ satisfies
\begin{equation}\label{discbound}
 D^*_N(z_1, \ldots, z_N) \leq 5.7 \sqrt{4.9 + \frac{\log \left((1-q)^{-1}\right)}{s}} \frac{\sqrt{s}}{\sqrt{N}}
\end{equation}
with probability at least $q$.
\end{theorem}

It is interesting that the quality of the discrepancy estimate in Theorem \ref{th1} \emph{improves} as the dimension $s$ increases; for example the necessary number $c(q,s)$ to have star-discrepancy bounded by $c(q,s) s^{1/2} N^{-1/2}$ with probability at least 90\% is 15.30 in dimension $s=1$, while it is only 12.65 in dimension $s=100$. However, neglecting this advantage of large dimensions in order to obtain a result which holds uniformly in $s$, one immediately obtains the following corollary.

\begin{cor} \label{co1}
For any $s \geq 1, ~N \geq 1$ and $q \in (0,1)$ a randomly generated $s$-dimensional point set $(z_1, \ldots, z_N)$ satisfies
\begin{equation}\label{discbound2}
 D^*_N(z_1, \ldots, z_N) \leq 5.7 \sqrt{4.9 + \log \left((1-q)^{-1}\right)} \frac{\sqrt{s}}{\sqrt{N}}
\end{equation}
with probability at least $q$.
\end{cor}

Theorem \ref{th1} shows that the probability that a random point set satisfies the discrepancy bound $c(q,s) s^{1/2} N^{-1/2}$ is extremely large already for moderate values of $c(q,s)$. The following table illustrates this fact, for $s=10$ and $s=100$.\\
\begin{center}
\begin{tabular}{|l||r|r|r|r|r|}
\hline
q & 0.01 & 0.5 & 0.9 & 0.99 & 0.999 \\
\hline
c(q,10) & 12.62 & 12.71 & 12.92 & 13.20 & 13.48 \\
\hline
c(q,100) & 12.62 & 12.63 & 12.65 & 12.68 & 12.71 \\
\hline
\end{tabular}
\end{center}
As the table shows, the probability that a random point set has ``small'' discrepancy in the sense that its discrepancy is bounded by $c s^{1/2} N^{-1/2}$ for some moderate $c$ (for example, $c=20$) is extremely large. This observation is an exciting counterpart of the fact that we do not have the slightest idea of how to construct point sets satisfying such discrepancy bounds, even for moderate $N$ and $s$. It should also be noted that calculating the star-discrepancy of a given (high-dimensional) point set is computationally very difficult, see \cite{gia, gnenp}. Hence, although our results show that the probability of a random point set having small discrepancy is very large, checking that a concrete point set satisfies such discrepancy bounds is in general (in high dimensions) a computationally intractable problem.

\section{Preliminaries}
Throughout the paper, $s \geq 1$ denotes the dimension and $\lambda$ denotes the $s$-dimensional Lebesgue measure. For $x,y \in [0,1]^s$, where $x = (x_1, \ldots, x_s)$ and $y = (y_1, \ldots, y_s)$, we write $x \leq y$ if $x_i \leq y_i, 1 \leq i \leq s$, and for any $x \in [0,1]^s$ we write $[0,x]$ for the set $\{y \in [0,1]^s:~0 \leq y \leq x\}$. Furthermore, we write $|A|$ for the number of elements of a set $A$.\\

The following Lemma \ref{gnew} of Gnewuch \cite[Theorem 1.15]{gb} is a central ingredient in the proof of our main result. For convenience we use the notation from \cite{gb} and \cite{gc}: For any $\delta \in (0,1]$ a set $\Gamma$ of points in $[0,1]^s$ is called a $\delta$-cover of $[0,1]^s$ if for every $y \in [0,1]^s$ there exist $x,z \in \Gamma \cup \{0\}$ such that $x \leq y \leq z$ and $\lambda ([0,z)) - \lambda ([0,x)) \leq \delta$. The number $\mathcal{N}(s,\delta)$ denotes the smallest possible cardinality of a $\delta$-cover of $[0,1]^s$.\\
Similarly, for any $\delta \in (0,1]$ a set $\Delta$ of pairs of points from $[0,1]^s$ is called a $\delta$-bracketing cover of $[0,1]^s$, if for every pair $(x,z) \in \Delta$ the estimate $\lambda ([0,z)) - \lambda ([0,x)) \leq \delta$ holds, and if for every $y \in [0,1]^s$ there exists a pair $(x,z)$ from $\Delta$ such that $x \leq y \leq z$. The number $\mathcal{N}_{[~]}(s,\delta)$ denotes the smallest possible cardinality of a $\delta$-bracketing cover of $[0,1]^s$.
\begin{lemma} \label{gnew}
For any $s \geq 1$ and $\delta \in (0,1]$
$$
\mathcal{N}(s,\delta) \leq (2e)^s (\delta^{-1} +1)^s
$$
and
$$
\mathcal{N}_{[~]}(s,\delta) \leq 2^{s-1} e^s (\delta^{-1} +1)^s.
$$
\end{lemma}
By Lemma \ref{gnew} for any $1 \leq k \leq K$ there exists a $2^{-k}$-cover of $[0,1]^s$, denoted by $\Gamma_k$, such that
\begin{equation*}
 |\Gamma_k| \leq (2e)^s (2^k +1)^s.
\end{equation*}
Furthermore we denote by $\Delta_K$ a $2^{-K}$- bracketing cover for which
\begin{equation*}
 |\Delta_K| \leq 2^{s-1} e^s (2^K +1)^s,
\end{equation*}
which also exists due to Lemma \ref{gnew}. Moreover we define $\Gamma_K$ as
\begin{equation*}
 \Gamma_K = \{v \in [0,1]^s:~(v,w) \in \Delta_K \textrm{~for some } w \}.
\end{equation*}
By definition for every $x \in [0,1]^s$ there exists a pair $(v_K,w_K)=(v_K(x),w_K(x))$ for which $(v_K,w_K) \in \Delta_K$ such that $v_K \leq x \leq w_K$ and
$$
\lambda ([0,w_K]) - \lambda ([0,v_K]) \leq \frac{1}{2^{K}}.
$$
Furthermore for every $k, ~2 \leq k \leq K$ and $\gamma \in \Gamma_k$ there exist $v_{k-1}=v_{k-1}(\gamma),w_{k-1}=w_{k-1}(\gamma),~v_{k-1},w_{k-1} \in \Gamma_{k-1} \cup \{0\}$, such that $v_{k-1} \leq \gamma \leq w_{k-1}$ and
$$
\lambda ([0,w_{k-1}]) - \lambda ([0,v_{k-1}]) \leq \frac{1}{2^{k-1}}.
$$
We define
\begin{eqnarray*}
p_K(x) & = & v_K(x)\\
p_{K-1}(x) & = & v_{K-1} (p_K(x)) = v_{K-1} (v_K(x))\\
p_{K-2}(x) & = & v_{K-2} (p_{K-1}(x)) = v_{K-2}(v_{K-1}(v_K(x))) \\
& \vdots& \\
p_1 (x) & = & v_1 (p_2(x)),
\end{eqnarray*}
and $p_{K+1}(x) = w_K(x), \qquad p_0(x) = 0$.
For $x,y \in [0,1]^s$ we set
$$
\overline{[x,y]} := \left\{ \begin{array}{ll} [0,y] \backslash [0,x] & \textrm{if $x \neq 0$,} \\ \textrm{}[0,y] & \textrm{if $x = 0,~y \neq 0$,} \\ \emptyset & \textrm{if $x=y=0$.} \end{array} \right.
$$
Then the sets $\overline{[p_k(x),p_{k+1}(x)]},\quad 1 \leq k \leq K,$ are disjoint, and we obtain
$$
\bigcup_{k=0}^{K-1} \overline{[p_k(x),p_{k+1}(x)]} \subset [0,x] \subset \bigcup_{k=0}^{K} \overline{[p_k(x),p_{k+1}(x)]}, \quad \forall x \in [0,1]^s.
$$
Hence for every $x,y \in [0,1]^s$
\begin{equation} \label{acc}
\sum_{k=0}^{K-1} \mathbf{1}_{\overline{[p_k(x),p_{k+1}(x)]}} (y) \leq \mathbf{1}_{[0,x]} (y) \leq \sum_{k=0}^{K} \mathbf{1}_{\overline{[p_k(x),p_{k+1}(x)]}} (y).
\end{equation}
Moreover, independent of $x$, we have for $0 \leq k \leq K$
\begin{eqnarray*}
\lambda \left(\overline{[p_k(x), p_{k+1}(x)]}\right) & \leq & \frac{1}{2^k}.
\end{eqnarray*}
For $0 \leq k \leq K$ we define $A_k$ to be the set of all sets of the form $\overline{[p_k(x),p_{k+1}(x)]},$ where $x \in [0,1]^s$. Then for $0 \leq k \leq K$, as a consequence of Lemma \ref{gnew}, we can bound the cardinality of $A_k$ by 
\begin{equation}\label{cover}
|A_k| \leq (2e)^s \left(2^{k+1} +1\right)^s.
\end{equation}
Note that all elements of $A_k$, where $0 \leq k \leq K$, have Lebesgue measure bounded by $2^{-k}$. This dyadic decomposition method was introduced in \cite{ady}, where it is described in more detail.\\

Let $X_1, \dots, X_N$ be independent, identically distributed (i.i.d.) random variables defined on some probability space $(\Omega, \mathcal{A}, \mathbb{P})$ having uniform distribution on $[0,1]^s$, and let $I \in A_k$ for some $k \geq 0$. Then the random variables $\mathbf{1}_I (X_1), \dots, \mathbf{1}_I (X_N)$ are i.i.d. random variables, having expected value $\lambda(I)$
and variance
\begin{equation} \label{var}
\lambda(I) - \lambda(I)^2 \leq \left\{ \begin{array}{ll} 2^{-k} (1 - 2^{-k}) & \textrm{for $k \geq 1$,} \\ 1/4 & \textrm{for $k=0$.} \end{array} \right.
\end{equation}
Since the $X_n$ are independent it follows that the random variable $\sum_{n=1}^N \mathbf{1}_I (X_n) $ has expected value $N \lambda(I)$ and variance $N (\lambda(I)-\lambda(I)^2)$.\\

In the proof of our main result we need two well-known results from probability theory, namely Bernstein's and Hoeffding's inequality. Bernstein's inequality states that for $Z_1, \dots, Z_N$ being i.i.d. random variables, satisfying $\E Z_n =0$ and $|Z_n| \leq C$ a.s. for some $C>0$,
$$
\mathbb{P} \left( \left| \sum_{n=1}^N Z_n \right| > t \right) \leq 2 \exp \left( - \frac{t^2}{2 \left( \sum_{n=1}^N \E Z_n^2 \right) + 2 C t/3} \right).
$$
By applying this inequality to the random variables $\mathbf{1}_I (X_n) - \lambda(I)$, we obtain 
$$
\mathbb{P} \left( \left| \sum_{n=1}^N \mathbf{1}_I (X_n) - N \lambda(I) \right| > t \right) \leq 2 \exp \left( - \frac{t^2}{2 \left( N \lambda (I) \left( 1 - \lambda(I) \right) \right) + 2 t/3} \right)
$$
for $t > 0$. Using (\ref{var}) we conclude
\begin{equation} \label{bernstein}
\mathbb{P} \left( \left| \sum_{n=1}^N \mathbf{1}_I (X_n) - N \lambda(I) \right| > t \right) \leq 2 \exp \left( - \frac{t^2}{2 N 2^{-k} (1-2^{-k}) + 2 t/3} \right) \quad \textrm{for}  \quad k \geq 2.
\end{equation}
For $k \in \{0,1\}$ we use Hoeffding's inequality, which yields
\begin{equation}\label{hoeffding}
\mathbb{P} \left( \left| \sum_{n=1}^N \mathbf{1}_I (X_n) - N \lambda(I) \right| > t \right) \leq 2 \exp \left( - \frac{2 t^2}{N} \right).
\end{equation}

\section{Proof of Theorem \ref{th1}}
Since the theorem is trivial for $N < 32 \left(s + \log \left((1-q)^{-1}\right) \right) < 5.7^2 \left(s + \log \left((1-q)^{-1}\right) \right)$ we assume that $N \geq 32 \left(s + \log \left((1-q)^{-1}\right) \right)$ and set
\begin{equation*}
 K = \left \lceil \frac{ \log_2 N  - \log_2\left(s + \log \left((1-q)^{-1}\right) \right)}{2} \right \rceil.
\end{equation*}
Then $K \geq 3$, and
\begin{equation} \label{Kest}
2^{-K} \in \left[ \frac{\sqrt{s + \log \left((1-q)^{-1}\right) }}{2 \sqrt{N}},\frac{\sqrt{s + \log \left((1-q)^{-1}\right)}}{\sqrt{N}} \right].
\end{equation}
Furthermore we have
\begin{equation}\label{sqrt}
\sqrt{sN} = N \frac{\sqrt{s}}{\sqrt{N}} \leq \frac{2^{-K+1} \sqrt{s}}{\sqrt{ s + \log \left((1-q)^{-1}\right)}} N.
\end{equation}
By choosing $t = c \sqrt{s N}$ for some $c > 0$, we conclude from \eqref{bernstein}, \eqref{hoeffding} and \eqref{sqrt} that for any $c>0$
\begin{eqnarray} 
& & \mathbb{P} \left( \left| \sum_{n=1}^N \mathbf{1}_I (X_n) - N \lambda(I) \right| > c \sqrt{s N} \right) \nonumber\\
& \leq & \left\{ \begin{array}{ll} 2 e^{- 2c^2 s}  & \textrm{for $k=0,1$}\\ 2 \exp \left( - \frac{c^2 s}{2 \cdot 2^{-k} (1-2^{-k}) + \frac{4 c 2^{-K} \sqrt{s}}{3 \sqrt{s + \log \left((1-q)^{-1}\right)}}} \right) & \textrm{for $2 \leq k \leq K$} . \end{array} \right.\label{pro}
\end{eqnarray}
Let $B_k, k = 0, \ldots, K$ be given as
\begin{equation}\label{bk}
B_k = \bigcup_{I \in A_k} \left( \left| \sum_{n=1}^N \mathbf{1}_I (X_n) - N \lambda(I) \right| > c_k \sqrt{s N} \right).
\end{equation}
The strategy of the proof is to find for any given $q \in (0,1)$ constants $c_k = c_k(q), k = 0, \ldots, K$ for which
\begin{equation*}
 \sum_{k = 0}^K \mathbb{P} \left( B_k \right) < 1 - q
\end{equation*}
holds for any given $q$.\\ 

First we consider the case $k = 0$. By \eqref{cover} we have that $|A_0| \leq (6 e)^s$. We choose
\begin{equation}
 c_0 = \sqrt{\frac{1 + \log(6)}{2} + \frac{\log\left( 8 (1-q)^{-1}\right)}{2 s}} \leq \frac{1}{\sqrt{2}} \sqrt{4.88 + \frac{\log \left((1-q)^{-1}\right)}{s}}, \label{c0approx}
\end{equation}
thus together with \eqref{pro} and \eqref{bk} it follows that
\begin{align*}
 \mathbb{P}(B_0) &\leq |A_0| 2 e^{- 2 c_0^2 s} \leq (6 e)^s 2 e^{- s (1 + \log(6))} \frac{(1-q)}{8} = \frac{1-q}{4}.
\end{align*}
Furthermore we get by \eqref{cover} that $|A_1| \leq (10 e)^s$ and with
\begin{equation}
 c_1 = \sqrt{\frac{1 + \log(10)}{2} + \frac{\log\left(8 (1-q)^{-1}\right)}{2 s}} \leq \frac{1}{\sqrt{2}} \sqrt{5.39 + \frac{\log\left((1-q)^{-1}\right)}{s}}, \label{c1approx}
\end{equation}
we obtain that
\begin{align*}
 \mathbb{P}(B_1) &\leq |A_1| 2 e^{- 2 c_1^2 s} \leq (10 e)^s 2 e^{- s (1 + \log(10))} \frac{(1-q)}{8}= \frac{1-q}{4}.
\end{align*}
Next we consider the case $2 \leq k \leq K$. By \eqref{pro} and \eqref{bk} we have
\begin{align}
 \mathbb{P}(B_k) &\leq |A_k| \cdot 2 \cdot \exp \left( - \frac{c_k^2 s}{2 \cdot 2^{-k} (1-2^{-k}) + \frac{4 c_k 2^{-K} \sqrt{s}}{3 \sqrt{s + \log \left((1-q)^{-1}\right)}}} \right).\label{ck}
\end{align}
We set
\begin{equation*}
 c_k = \sqrt{1 + \log(2 (2^{k+1} + 1)) + \frac{\log\left(2^{(k + 1)} (1-q)^{-1} \right)}{s} } \sqrt{2 \cdot 2^{-k} (1 - 2^{-k}) + \frac{2.08 \cdot 4 \cdot 2^{-K}}{3}}
\end{equation*}
hence we get that
\begin{align*}
&\left| \frac{c_k \sqrt{s}}{\sqrt{s + \log \left((1-q)^{-1}\right)}} \right|\\ 
\leq &\left| \frac{\sqrt{1 + \log(2 (2^{k+1} + 1)) + \log(2^{(k + 1)}) + \frac{\log\left( (1-q)^{-1}\right)}{s} } \sqrt{2 \cdot 2^{-k} (1 - 2^{-k}) + \frac{2.08 \cdot 4 \cdot 2^{-K}}{3}}}{\sqrt{1 + \frac{\log\left((1 - q)^{-1}\right)}{s}}} \right|\\
 \leq &\left| \sqrt{1 + \log(2 (2^{k+1} + 1)) + \log(2^{(k + 1)})} \sqrt{2 \cdot 2^{-k} (1 - 2^{-k}) + \frac{2.08 \cdot 4 \cdot 2^{-K}}{3}}\right| \\
 \leq &2.08
\end{align*}
for $2 \leq k \leq K$. Thus by \eqref{ck} we obtain
\begin{align*}
 \mathbb{P}(B_k) &\leq |A_k| \cdot 2 \cdot \exp \left( - \frac{c_k^2 s}{2 \cdot 2^{-k} (1-2^{-k}) + \frac{4 c_k 2^{-K} \sqrt{s}}{3 \sqrt{s + \log\left((1 - q)^{-1}\right)}}} \right)\\
 &\leq  (2 e)^s (2^{(k + 1)} + 1)^s \cdot 2 \cdot \exp \left( -s \left( 1 + \log(2 (2^{(k + 1)} + 1)) \right) \right) \frac{1 - q}{2^{(k + 1)}}\\
 &= \frac{1 - q}{2^{k}}.
\end{align*}
Summing up the estimated probabilities gives
\begin{equation*}
 \sum_{k = 0}^K \mathbb{P}(B_k) \leq \left( \frac{3}{4} + \sum_{k = 3}^K 2^{-k} \right) (1 - q) < 1 - q.
\end{equation*}
Therefore with at least probability $q$, a realization $X_1(\omega), \ldots, X_n(\omega)$ is such that
\begin{equation*}
 \omega \notin \bigcup_{k = 0}^K B_k.
\end{equation*}
We denote by $z_n$ a point set which is defined by such a realization, i.e.\
\begin{equation*}
 z_n = X_n(\omega), \quad 1 \leq n \leq N, \qquad \textrm{for some~}\omega \notin \bigcup_{k = 0}^K B_k.
\end{equation*}
Set $\lambda_k = \sqrt{2 \cdot 2^{-k} (1 - 2^{-k}) + 2.08 \cdot 4 \cdot 2^{-K} / 3}$. Then
\begin{align}
 \sum_{k = 2}^{K} c_k &= \sum_{k = 2}^{K} \lambda_k \sqrt{1 + \log(2 (2^{k+1} + 1)) + \frac{\log(2^{(k + 1)} (1-q)^{-1}  )}{s} } \notag \\
 &\leq \sum_{k = 2}^{K} \lambda_k \sqrt{1 + \log 2 +  \log(2^{k+1}) + 0.12 + \frac{\log(2^{(k + 1)})}{s} + \frac{\log\left((1-q)^{-1}\right)}{s} } \notag \\
&\leq \sum_{k = 2}^{K} \lambda_k \sqrt{1.12 + \left( 2k + 3 \right) \log 2 + \frac{\log\left((1-q)^{-1}\right)}{s} }  \notag \\
&\leq \sqrt{1.12 + 7 \log 2 + \frac{\log\left((1-q)^{-1}\right)}{s}} \frac{1}{\sqrt{2}} \sum_{k = 2}^{K} \sqrt{k} \lambda_k \notag\\
&\leq  3.28 \sqrt{5.98 + \frac{\log\left((1-q)^{-1}\right)}{s} }. \label{ckapprox}
\end{align}
Therefore we obtain by using \eqref{c0approx}, \eqref{c1approx} and \eqref{ckapprox}
\begin{align}
 \sum_{k = 0}^K c_k \leq &\frac{1}{\sqrt{2}} \sqrt{4.88 + \frac{\log\left((1-q)^{-1}\right)}{s}} + \frac{1}{\sqrt{2}} \sqrt{5.39 + \frac{\log\left((1-q)^{-1}\right)}{s}}\notag\\ 
& + 3.28 \sqrt{5.98 + \frac{\log\left((1-q)^{-1}\right)}{s}} \label{sumapprox}
\end{align}
Applying \eqref{acc}, \eqref{Kest}, \eqref{sumapprox} and Jensen's inequality we obtain
\begin{align*}
 \sum_{n = 1}^N \mathbf{1}_{[0,x]}(z_n) &\leq \sum_{k = 0}^{K} \sum_{n = 1}^N \mathbf{1}_{\overline{[p_k(x),p_{k+1}(x)]}} (z_n)\nonumber\\
& \leq N \lambda([0,w_K(x)]) + \sqrt{sN} \sum_{k=0}^K c_k \nonumber\\
& \leq  N \lambda([0,x]) + N \lambda(\overline{[x,w_K(x)]}) +  \sqrt{sN} \sum_{k=0}^K c_k\nonumber\\
& \leq N \lambda([0,x]) + N \frac{\sqrt{s + \log \left((1-q)^{-1}\right)}}{\sqrt{N}} + \sqrt{sN} \sum_{k=0}^K c_k\nonumber\\
& \leq N \lambda([0,x]) + 5.7 \sqrt{4.9 + \frac{\log\left((1-q)^{-1}\right)}{s}} \sqrt{sN}.
\end{align*}
Similarly a lower bound is given by
\begin{align*}
\sum_{n=1}^N \mathbf{1}_{[0,x]} (z_n) & \geq N \lambda([0,p_K(x)]) - \sqrt{sN} \sum_{k=0}^{K-1} c_k \\
& \geq N \lambda([0,x]) - 5.7 \sqrt{4.9 + \frac{\log\left((1-q)^{-1}\right)}{s} } \sqrt{sN}.
\end{align*}
Combining the above bounds we finally arrive at
\begin{equation*}
 D_N^* (z_1, \ldots, z_n) \leq 5.7 \sqrt{4.9 + \frac{\log\left((1-q)^{-1}\right)}{s}} \frac{\sqrt{s}}{\sqrt{N}}.
\end{equation*}

\end{document}